\documentclass[12pt,letterpaper,reqno]{amsart}

\usepackage{times} \usepackage[T1]{fontenc} \usepackage{mathrsfs} \usepackage{latexsym} \usepackage[dvips]{graphics}
\usepackage{epsfig}
\usepackage{amsmath,amsfonts,amsthm,amssymb,amscd} \input amssym.def \input amssym.tex

 \newcommand\be{\begin{equation}} \newcommand\ee{\end{equation}}
\newcommand\bea{\begin{eqnarray}} \newcommand\eea{\end{eqnarray}} \newcommand\bi{\begin{itemize}}
\newcommand\ei{\end{itemize}} \newcommand\ben{\begin{enumerate}} \newcommand\een{\end{enumerate}}
\newcommand\bc{\begin{center}} \newcommand\ec{\end{center}} \newcommand\ba{\begin{array}} \newcommand\ea{\end{array}}

\newtheorem{thm}{Theorem}[section] \newtheorem{conj}[thm]{Conjecture} 
 \newtheorem{prop}[thm]{Proposition} \newtheorem{exa}[thm]{Example}
\newtheorem{defi}[thm]{Definition}

\theoremstyle{definition} \newtheorem{rek}[thm]{Remark}

\begin{document}

\title[A Cauchy-Davenport type result for arbitrary regular graphs] {A Cauchy-Davenport type result for arbitrary regular graphs}

\author{Peter Hegarty} \address{Mathematical Sciences, Chalmers University Of Technology and University of Gothenburg,
41296 Gothenburg, Sweden} \email{hegarty@chalmers.se}


\date{\today}

\begin{abstract} Motivated by the Cauchy-Davenport theorem for sumsets, and its interpretation in terms of
Cayley graphs, we prove the following main result : There is a universal constant $\epsilon > 0$ such that, if
$\mathscr{G}$ is a connected, regular graph on $n$ vertices, then either every pair of vertices can be
connected by a path of length at most three, or the number of pairs of such vertices is at least $1+\epsilon$ times the number of edges in $\mathscr{G}$. We discuss a range of further questions to which this result gives rise.
\end{abstract}


\maketitle

\setcounter{equation}{0}

\setcounter{equation}{0}

\section{Introduction and Statement of Results}

Let $A$ be a subset of an abelian group $G$, written additively, and $h$ a positive integer. The $h$-fold sumset $hA$ is defined as
\be
hA = \{g \in G : g = a_1 + \cdots + a_h \; {\hbox{for some $a_1,...,a_h \in A$}} \}.
\ee
We say that $A$ is a {\em basis} for $G$ if $hA = G$ for some $h$. The cardinality of a set $S$ will be denoted $|S|$.
The following is a (special case of a) fundamental result in the theory of sumsets :

\begin{thm} ({\bf Cauchy-Davenport}) Let $p$ be a prime and $A$ a subset of
$\mathbb{Z}_{p}$. Then
\be
|hA| \geq \min\{p, h|A| -(h-1)\}.
\ee
\end{thm}

There is a well-known generalisation of this result to arbitrary abelian groups, due to Kneser \cite{K}, but that is not what is of primary interest to us here. Instead, we
are interested in interpreting the Cauchy-Davenport result in terms of graphs. First,
recall the following definition :

\begin{defi}
Let $G$ be an abelian group and $S$ a subset of $G$. The {\em Cayley graph}{\footnote{Throughout this paper, the letter $G$ will be reserved to denote an abelian group and graphs will be denoted by the scripted letter $\mathscr{G}$.}}
$\mathscr{G} = \mathscr{G}(G,S)$ is the directed graph whose vertices are the elements of $G$ and whose edges consist of the ordered pairs $(g_1,g_2)$ such that $g_2 - g_1 \in S$.
\end{defi}

Note that the graph $\mathscr{G}(G,S)$ is strongly connected if and only if $S$ is a basis for $G$.
For simplicity, let us assume this, plus that the set $S$ is symmetric, i.e.: $S = -S$, and contains the identity element of $G$. Then
we can think of the Cayley graph as being undirected, with a loop at every vertex. In this case, let
$\mathscr{G}_{0}(G,S)$ be the part of $\mathscr{G}(G,S)$ with all the loops removed. For the rest of the paper, we shall only deal with undirected, loopless graphs. Now consider the following definition :

\begin{defi}
Let $\mathscr{G}$ be a graph on vertex set $V$ and $h$ a positive integer. We denote by $h\mathscr{G}$ the graph on vertex set $V$ such that $\{v_1,v_2\}$ is an edge in
$h\mathscr{G}$ if and only there is a path joining $v_1$ to $v_2$ in $\mathscr{G}$ of
length at most $h$. We shall call $h\mathscr{G}$ the {\em $h$-fold sumgraph}{\footnote{This definition is well-known in the literature, though it is standard to use multiplicative notation instead, which is natural when one thinks in terms of the adjacency matrix of the graph. So what we are calling the $h$-fold sumgraph
$h\mathscr{G}$ is usually referred to as the {\em $h$-th power of $\mathscr{G}$} and denoted $\mathscr{G}^h$. Observe that if we
add a loop at each vertex of $\mathscr{G}$ and let
$A$ be the adjacency matrix of the resulting graph, then $|\mathscr{E}(\mathscr{G}^h)|$ is just half the number of non-zero off-diagonal entries in
$A^h$.
\par For the remainder of this note we shall retain our additive notation and terminology for graphs so as to emphasise
the connection to sumsets.}} of
$\mathscr{G}$ and denote its set of edges by $h\mathscr{E} = \mathscr{E}(h\mathscr{G})$.
\end{defi}

Then the following is an immediate consequence of the Cauchy-Davenport \\ theorem :

\begin{thm}
Let $p$ be a prime, $A$ a subset of $\mathbb{Z}_{p}$ such that $0 \in A$ and $A = -A$.
Let $\mathscr{G} = \mathscr{G}_{0}(\mathbb{Z}_{p},A)$ be the Cayley graph of $A$, minus all loops. Then for
every positive integer $h$,
\be
|h\mathscr{E}| \geq \min \left\{ \left( \begin{array}{c} p \\ 2 \end{array} \right), h|\mathscr{E}| \right\}.
\ee
\end{thm}

The question which motivated this paper is whether anything like this result is true for more general connected graphs. More precisely, the feature of Theorem 1.4 that we are interested in generalising is the fact that the (edge) sizes of the graphs $h\mathscr{G}$ grow at least linearly in $h$, as long as $\mathscr{G}$ isn't already too dense. As we shall show below, it is hopelessly optimistic to hope for anything like this phenomenon in arbitrary connected graphs. However, Cayley graphs have the very important property that they are regular. Our main result is the following partial generalisation of Theorem 1.4 :

\begin{thm}
There is a universal constant $\epsilon > 0$ such that if $\mathscr{G}$ is a regular, connected graph on $n$ vertices, then
\be
|3\mathscr{E}| \geq \min \left\{ \left( \begin{array}{c} n \\ 2 \end{array} \right),  (1+\epsilon)|\mathscr{E}| \right\}.
\ee
In fact, we can take $\epsilon$ to be the unique positive root of the equation
\be
\epsilon = \frac{1}{4} (1-\sqrt{\epsilon})^3,
\ee
i.e.: $\epsilon \approx 0.087...$.
\end{thm}

We were surprised by the simplicity and elegance of this result, which is why we
considered it worth mentioning. Of course, it is unsatisfactory in many respects so some
detailed remarks are in order :
\\
\\
1. The obvious problem with our result is that it cannot be used recursively to obtain
estimates for the growth of $h$-fold sumgraphs for arbitrary $h$. This is because, even
if the graph $\mathscr{G}$ is regular, then the graphs $h\mathscr{G}$ need not be, for any $h > 1$ (note that regularity is preserved for Cayley sum-graphs). Thus it remains to obtain a generalisation of Theorem 1.5 to
$h$-fold sumgraphs for arbitrary $h$. Note that, for fixed degree, the sumgraphs
$h\mathscr{G}$ grow at least linearly $\lq$on average' until the graph becomes complete. This is a trivial observation, but a more precise result is contained in the next proposition.
Recall that the {\em diameter} of a graph is the smallest $\delta > 0$ such that any pair of vertices are connected by a path of
length at most $\delta$. In the notation of Definition 1.3, the diameter of
a graph $\mathscr{G}$ on $n$ vertices is the smallest $h$ such that $h\mathscr{G} =
K_{n}$, the complete graph. Now we have

\begin{prop}
Let $\mathscr{G}$ be a connected graph on $n$ vertices and of minimal degree $d$. Then
\be
{\hbox{diam}}(\mathscr{G}) \leq \frac{3n-(d+3)}{d+1}.
\ee
\end{prop}

2. However, the growth of sumgraphs can certainly be irregular. In particular, and this is the most natural thing to ask about, there is no constant $\epsilon^{\prime} > 0$ such that the analogue of Theorem 1.5 holds for
$2$-fold sumgraphs. To see this consider the following \\ example :

\begin{exa}
Fix $d > 0$ and let $n$ be a multiple of $d+1$, say $n = m(d+1)$. Let $\mathscr{G} =
\mathscr{G}_{d,m}$ be the following graph on $n$ vertices : Partition the vertex set
$V$ into $m$ disjoint subsets of size $d+1$, say $V_1,V_2,...,V_m$. For each
$i = 1,...,m$ pick two vertices $v_{i1}, v_{i2} \in V_i$. Now the graph
$\mathscr{G}_{d,m}$ contains the following edges :
\par (i) for each $i = 1,...,m$, insert all edges among the vertices of $V_{i}$,
except the edge $\{v_{i1},v_{i2}\}$.
\par (ii) for each $i = 1,...m-1$, insert the edge $\{v_{i1},v_{(i+1),2}\}$, and then
finally add the edge $\{v_{m1},v_{12}\}$.
\\
\\
Clearly, this graph is connected and $d$-regular, so
\be
|\mathscr{E}_{d,m}| = \left( \frac{d}{2} \right) n.
\ee
However, one easily checks that
\be
|2\mathscr{E}_{d,m} \backslash \mathscr{E}_{d,m}| = \left( \frac{2d-1}{d+1} \right) n,
\ee
so that $|2\mathscr{E}_{d,m}| = (1+o_{d}(1))|\mathscr{E}_{d,m}|$. Note also that for this
graph one may check that
\be
|3\mathscr{E}_{d,m} \backslash \mathscr{E}_{d,m}| = \left( \frac{d^{2}+4}{d+1} \right) n,
\ee
so that $|3\mathscr{E}_{d,m}| = (3 - o_{d}(1))|\mathscr{E}_{d,m}|$.
\end{exa}

Considering this example naturally leads one to asking for more precise extremal results. We
believe that the graphs $\mathscr{G}_{d,m}$ are essentially extremal for $2$-fold sumgraphs, but these latter objects are still somewhat mysterious to us. Motivated by (1.8), we are prepared at this stage to conjecture the following :

\begin{conj}
Let $d,n$ be positive integers. If $\mathscr{G}$ is a
$d$-regular graph on $n$ vertices, then either $2\mathscr{G} = K_n$ or
\be
|2\mathscr{E} \backslash \mathscr{E}| \geq \left( 2 - o_{d}(1) \right) n.
\ee
\end{conj}

Note that, in the notation of this conjecture, if $n \geq d+2$ then trivially $|2\mathscr{E}\backslash \mathscr{E}|
\geq n/2$, since every vertex is connected to at least one non-neighbor by a path of length two. Hence there
is a factor of four separating (asymptotically) 
the trivial lower bound for $|2\mathscr{E}\backslash \mathscr{E}|$ and what we
conjecture to be the truth.
\par Neither is it
clear to us whether the graphs $\mathscr{G}_{d,m}$ are essentially extremal for $3$-fold sumgraphs. The question here
is what is the best-possible choice of the constant $\epsilon$ in Theorem 1.5 ?
By (1.9), we can't take $\epsilon > 2$. Indeed, the same conclusion could be drawn by considering the Cayley graph of an arithmetic progression.
\par Also, note that the graphs $\mathscr{G}_{d,m}$ are certainly not close to being extremal
sumgraphs in general. This is because it is easy to see that $\mathscr{G}_{d,m}$ has
diameter $m+1 = \frac{n}{d+1} + 1$, whereas from the proof of Proposition 1.6 we will
easily be able to construct examples which show that the upper bound in
(1.6) is essentially best-possible, even for regular graphs (see Remark 3.1). Hence, we suspect that the extremal problem for sumgraphs in general might be quite hard.
\\
\\
3. Finally, note that there doesn't seem to be any hope of obtaining meaningful generalisations of our results to graphs which are not regular. For example, let $n$ be a positive integer and let $\mathscr{G}_n$ be the graph on $n$ vertices
which is the union of a complete subgraph on $\lfloor n^{3/4} \rfloor$ vertices and a path of length
$n - \lfloor n^{3/4} \rfloor$ which is joined to the complete subgraph at one vertex. This graph is connected and contains $\Theta (n^{3/2})$ edges but, for any fixed $h$, the $h$-fold sumgraph contains only $\Theta_{h}(n)$ additional edges.
\\
\\
The rest of the paper is organised as follows. Sections 2 and 3 are devoted to the proofs and discussion of
Theorem 1.5 and Proposition 1.6 respectively. Section 4 contains a quick recap of unresolved issues.

\setcounter{equation}{0}

\section{Proof of Theorem 1.5}

{\bf Notation.} If $\mathscr{G}$ is a graph and $X \subseteq V(\mathscr{G})$, then $N(X)$ will denote the set
of all neighbors of the vertices in $X$. If $X$ is a singleton set, say $X = \{x\}$, then we simply write
$N(x)$.
\\
\\
Let $d,n$ be positive integers and let $\mathscr{G}$ be a connected, $d$-regular graph on $n$ vertices. Let $\epsilon > 0$ be the solution of (1.5) and suppose that
$|3\mathscr{G}| < (1+\epsilon)|\mathscr{G}|$. We must show that $3\mathscr{G} = K_n$. Since $2\mathscr{E} \subseteq
3\mathscr{E}$, we can first of all deduce that $|2\mathscr{E}| < (1+\epsilon)|\mathscr{E}|$. We present the argument
in a sequence of steps.
\\
\\
{\sc Step 1} : Set $\epsilon_{1} := \sqrt{\epsilon}$. For each $v \in V(\mathscr{G})$, let
\be
T_{v} := \{w \in V(\mathscr{G}) : \{v,w\} \in 2\mathscr{E}\backslash \mathscr{E} \}
\ee
and let
\be
V_1 := \{v \in V(\mathscr{G}) : |T_{v}| < \epsilon_{1} d. \}
\ee
Since, by assumption,
\be
\epsilon d n > 2 \times |2\mathscr{G} \backslash \mathscr{G}| = \sum_{v} |T_{v}|,
\ee
it follows easily that
\be
|V_{1}| > (1-\epsilon_{1}) n.
\ee
{\sc Step 2} : Let $v \in V_1$. Set $A_v := N(v)$, $B_v := N(N(v))$ and $C_v := B_v \backslash (\{v\} \cup A_v)$. If the set $C_v$ were empty then, since the graph is connected, it would imply that $V(\mathscr{G}) = \{v\} \cup A_v$ and
hence that $2\mathscr{G} = K_n$. So we may assume that $C_v$ is non-empty. If $c \in C_v$ then there is a path $v \rightarrow a
\rightarrow c$ in $\mathscr{G}$, for some $a \in A_v$, hence $\{v,c\} \in 2\mathscr{E}$. By definition of the set
$V_1$, it follows that
\be
|C_v| < \epsilon_1 d.
\ee
Set $D_v := V(\mathscr{G}) \backslash (\{v\} \cup A_v \cup C_v)$. Suppose this set were empty. Since $\mathscr{G}$ is connected, it would imply that $V(\mathscr{G}) = \{v\} \cup A_v \cup C_v$. We claim that, in this case, $3\mathscr{G} = K_n$. We need to show that any two vertices in $\{v\} \cup A_v \cup C_v$ can be connected by a path of length at most 3. This is obvious unless both vertices lie in $C_v$. Consider a pair of such vertices, say $c_1$ and $c_2$. Our
assumptions say that $N(c_1) \cup N(c_2) \subseteq A_v \cup C_v$. But by (2.5), $d$-regularity and the fact that $\epsilon_1 < 1/2$ (see the statement of Theorem 1.5), it follows that $c_1$ and $c_2$ must have a common neighbor in $A_v$. Hence $\{c_1,c_2\} \in
2\mathscr{E}$, in fact.
\par So we may now assume that the set $D_v$ is non-empty. There must be at least one edge between $C_v$ and $D_v$. For any such edge, say $\{c_v,d_v\}$, we know by (2.5) that
at least $(1-\epsilon_1)d$ of the neighbors of $c_v$ lie in $A_v \cup D_v$. Let $\mathscr{C}_v$ be the set of
vertices in $C_v$ with at least one neighbor in $D_v$ and set
\be
\alpha_v := \frac{1-\epsilon_1}{d} \times \max \{|N(c_v) \cap A_v| : c_v \in \mathscr{C}_v \}.
\ee
In the steps to come, we consider the following two cases, at least one of which must obviously apply :
\\
\\
{\em Case 1} : For at least half of all $v \in V_1$, one has $\alpha_v \leq \frac{1}{2}$.
\\
{\em Case 2} : For at least half of all $v \in V_1$, one has $\alpha_v > \frac{1}{2}$.
\\
\\
{\sc Step 3} : Suppose {\em Case 1} holds. Let $V_2 := \{v \in V_1 : \alpha_v \leq \frac{1}{2}\}$. For each $v \in V_2$,
pick any vertex $c_v \in \mathscr{C}_v$. Then there are at least $\frac{1}{2}(1-\epsilon_1)d$ choices for an
edge $\{c_v,d_v\}$ such that $d_v \in D_v$. Notice that,
for any choice of $d_v$, there is a path $v \rightarrow a_v \rightarrow c_v \rightarrow d_v$ in $\mathscr{G}$, for some $a_v \in A_v$. Hence
$\{v,d_v\} \in 3\mathscr{E} \backslash \mathscr{E}$. Summing over all $v \in V_2$ and noting that any given pair of vertices is counted at
most twice, it follows that
\be
|3\mathscr{E} \backslash \mathscr{E}| \geq \frac{1}{2} \cdot |V_2| \cdot \frac{1}{2}(1-\epsilon_1) d \geq
\frac{1}{2} (1-\epsilon_1)^2 |\mathscr{E}| > \epsilon |\mathscr{E}|,
\ee
contradicting our assumptions.
\\
\\
{\sc Step 4} : Suppose {\em Case 2} holds. Let $V_3 := \{v \in V_1 : \alpha_v > \frac{1}{2} \}$, so that
$|V_3| \geq \frac{1}{2}|V|$. Let $v \in V_3$ and fix a choice of a vertex $c_v \in \mathscr{C}_v$ such that $c_v$ has at least
$\frac{1}{2}(1-\epsilon_1)d$ neighbors inside $A_v$. Let $d_v$ be any neighbor of $c_v$ inside $D_v$. Observe that all the neighbors of $d_v$ lie inside $C_v \cup D_v$.
Hence, by (2.5), there are at least $(1-\epsilon_1)d$ choices for a vertex $e_v \in N(d_v) \cap D_v$. For any
such vertex $e_v$ and any vertex $a_v \in N(c_v) \cap A_v$, there is a path in the graph $a_v \rightarrow
c_v \rightarrow d_v \rightarrow e_v$. Hence $\{a_v,e_v\} \in 3\mathscr{E} \backslash \mathscr{E}$.
Therefore, if we set
\be
\mathbb{S} := \sum_{v \in V_3} \# \{ \{a_v,e_v\} \in 3\mathscr{G} \backslash \mathscr{G} :
a_v \in A_v, \; e_v \in D_v \},
\ee
then we have
\be
|\mathbb{S}| \geq |V_3| \times \frac{(1-\epsilon_1)d}{2} \times (1-\epsilon_1)d \geq \left( \frac{(1-\epsilon_1)^{3}}{2} \right) \frac{d^2 n}{2}.
\ee
On the other hand, since
since $v \in N(a_v)$ always, any pair of vertices can appear in the sum at most $2d$ times. It follows that
\be
|3\mathscr{E} \backslash \mathscr{E}| \geq  \frac{1}{4} (1-\epsilon_1)^{3} |\mathscr{E}| \geq \epsilon |\mathscr{E}|,
\ee
which again contradicts our assumptions, and completes the proof of the theorem.

\setcounter{equation}{0}

\section{Proof of Proposition 1.6}

Let $\mathscr{G}$ be a connected graph on $n$ vertices of minimal degree $d$. Let $\delta$ be the diameter of $\mathscr{G}$ and let
$v,w$ be a pair of vertices such that a shortest path between them has length exactly $\delta$. Let such a path be
\be
v_1 = v \rightarrow v_2 \rightarrow \cdots \rightarrow v_{\delta} \rightarrow v_{\delta+1} = w.
\ee
Let $A$ be the set of vertices along the path and $B := V(\mathscr{G}) \backslash A$. Using the fact that there is no shorter path
in $\mathscr{G}$ between $v$ and $w$, we shall count in two ways the
number $e(A,B)$ of edges in $\mathscr{G}$ between $A$ and $B$. On the one hand, this fact implies that there are no edges between the vertices along the path other than those in the
path itself. Since $\mathscr{G}$ has minimal degree $d$, it follows that
\be
e(A,B) \geq (\delta - 1)\cdot (d-2) + 2 \cdot (d-1) = (\delta + 1)(d-2) + 2.
\ee
On the other hand, the absence of a shorter path between $v$ and $w$ means that no vertex in $B$ can be joined to more than three vertices of $A$ (and if it joined to exactly three of them, then they must be adjacent along the path (3.1)).
Hence,
\be
e(A,B) \leq 3|B| = 3(n-\delta-1).
\ee
From (3.2) and (3.3) one easily deduces (1.6).

\begin{rek}
The proof just given can be easily adapted to construct explicit examples of graphs which show that the upper bound in (1.6) is essentially best-possible. Let $d \geq 5$ be odd for simplicity and choose a non-negative integer $k$. Let
\be
a := 3(k+1), \;\;\; b := k(d-2) + 2(d-1) = (k+2)(d-2) + 2, \;\;\; n := a+b.
\ee
We construct a $d$-regular graph on $n$ vertices as follows. The vertices of $\mathscr{G}$ are partitioned into
two disjoint sets $A$ and $B$ such that $|A| = a$, $|B|=b$. Denote
\be
A := \{v_1,...,v_a\}, \;\;\; B:= \{w_1,...,w_b\}.
\ee
The graph $\mathscr{G}$ will contain the following edges :
\\
\\
{\sc Type 1} : The edges of the path $v_1 \rightarrow v_2 \rightarrow \cdots \rightarrow v_a$.
\\
{\sc Type II} : All edges $\{v_i,w_j\}$ such that $1 \leq i \leq 3$ and $1 \leq j \leq d-1$, except the edges
$\{v_2,w_{d-1}\}$ and $\{v_3,w_1\}$.
\\
{\sc Type III} : All edges $\{v_{(a+1)-i}, w_{(b+1)-j}\}$ such that $1 \leq i \leq 3$ and $1 \leq j \leq d-1$,
except the edges $\{v_{a-1},w_{b+2-d}\}$ and $\{v_{a-2},w_{b}\}$.
\\
{\sc Type IV} : All edges $\{v_{3r+s},w_{(d-2)r+1+t}\}$ such that $1 \leq r \leq k$, $1 \leq s \leq 3$ and
$1 \leq t \leq d-2$.
\\
{\sc Type V} : The complete subgraph on the vertices $w_1,...,w_{d-1}$, minus a perfect matching on
the $d-3$ vertices $w_2,...,w_{d-2}$.
\\
{\sc Type VI} : The complete subgraph on the vertices $w_{b+2-d},...,w_{b}$, minus a perfect matching on the
$d-3$ vertices $w_{b+3-d},...,w_{b-1}$.
\\
{\sc Type VII} : For each $1 \leq r \leq k$, the complete subgraph on the vertices \\
$w_{r(d-2)+2},...,w_{(r+1)(d-2)+1}$.
\\
\\
One can readily check that this graph is indeed $d$-regular and, for $k \gg 0$, of diameter $a-1 := \delta$. Moreover,
for $d \geq 5$ one has
\be
\delta = \lfloor \frac{3n-(d+3)}{d+1} \rfloor.
\ee
\end{rek}

\vspace*{1cm}




\setcounter{equation}{0}

\section{Concluding Remarks}

There are two obvious directions in which the results of this paper need to be improved upon. The first is to generalise them to $h$-fold sumgraphs for arbitrary $h$, and in particular to understand better the most natural case when $h=2$. The second is to sharpen them, in particular to obtain the best-possible constant
$\epsilon$ in Theorem 1.5. Both directions naturally lead in turn to
Freiman-type inverse problems, where one wishes to say something about the $\lq$structure' of
regular, connected graphs whose sumgraphs grow slowly.

\setcounter{equation}{0}

\section*{Acknowledgements}
This work was performed while the author was visiting the Institute for Pure and Applied Mathematics (IPAM) at UCLA, and I thank them for their hospitality. My research is partly supported by a grant from the Swedish Science Research Council (Vetenskapsr\aa det).

\end{document}